\documentclass[english]{amsart}
\usepackage[T1]{fontenc}
\usepackage[latin1]{inputenc}
\usepackage{a4wide}
\usepackage{babel}
\usepackage{graphics}
\usepackage{psfrag}
\usepackage{setspace}
\makeatletter

\providecommand{\LyX}{L\kern-.1667em\lower.25em\hbox{Y}\kern-.125emX\@}

 \theoremstyle{plain}    
 \newtheorem{thm}{Theorem}[section]
 \numberwithin{equation}{section} 
 \numberwithin{figure}{section} 
 \theoremstyle{plain}    
 \newtheorem{lem}[thm]{Lemma} 

\def\H{\text{!`}}
\def\HH{\text{\em !`}}
\def\la{\lambda}

\makeatother
\begin{document}

\title{Proof of two conjectures of Zuber on fully packed loop configurations}

\author[Fabrizio Caselli]{Fabrizio Caselli$^\dagger$}

\author[Christian Krattenthaler]{Christian Krattenthaler$^\dagger$}

\address{Institut Girard Desargues,
Universit\'e Claude Bernard Lyon-I,
21, avenue Claude Bernard,
F-69622 Villeurbanne Cedex, France}
\thanks{$^\dagger$Research supported  
by EC's IHRP Programme, grant HPRN-CT-2001-00272, 
``Algebraic Combinatorics in Europe"}



\begin{abstract}
Two conjectures of Zuber [``On the counting of fully packed loops 
configurations.
Some new conjectures,'' preprint] on the enumeration of configurations in the
fully packed loop model on the square grid with periodic boundary
conditions, which have a prescribed linkage pattern, are proved.
Following an idea of de Gier [``Loops, matchings and alternating-sign 
matrices,'' \emph{Discrete Math.}, to appear],
the proofs are based on bijections between such fully packed loop 
configurations and rhombus tilings, and the hook-content formula
for semistandard tableaux.
\end{abstract}

\maketitle

\section{Introduction\label{intro}}

\begin{figure}[h]

{\centering \includegraphics{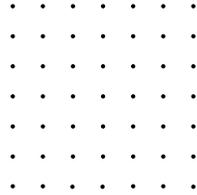} \par}

\caption{The square grid \protect\( Q_{7}\protect \)\label{fpl22}}
\end{figure}

The {\it fully packed loop model\/} (FPL model, for short;
see for example \cite{FPL})
is a model of (not necessarily closed) 
polygons on a lattice such that each vertex of the lattice
is on exactly one polygon. Whether or not these polygons are closed,
they will be also referred to as {\it loops}. 
Throughout this article, we consider this model on the 
square grid of side length \( n-1 \), which we denote by $Q_n$. See
Figure~\ref{fpl22} for a picture of $Q_7$.
\medskip

\medskip

\begin{figure}[h]
{\centering \includegraphics{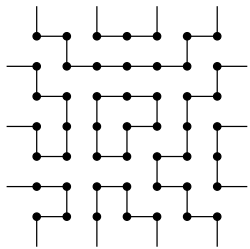} \par}

\caption{An FPL configuration on $Q_7$ with periodic boundary 
conditions\label{fpl24}}
\end{figure}

\begin{figure}
\psfrag{1}{$1$}
\psfrag{2}{$2$}
\psfrag{3}{$3$}
\psfrag{4n-1}{\kern-16pt$4n-1$}
\psfrag{4n}{$4n$}
\psfrag{...}{\dots}

{\centering \includegraphics{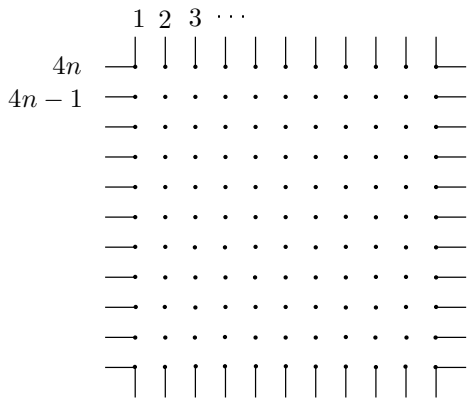} \par}

\caption{The labelling of the external links\label{fpl1}}
\end{figure}

\noindent
The polygons consist of horizontal or vertical edges connecting vertices
of $Q_n$, and
edges that lead outside of $Q_n$ from a vertex along the border of
$Q_n$, see Figure~\ref{fpl24} for an example of an allowed
configuration in the FPL model. 
We call the edges that stick outside of $Q_n$ {\it external
links}. The reader is referred to Figure~\ref{fpl1} for an
illustration of the external links of the square $Q_{11}$. (The labels
should be ignored at this point.) It should be noted that the four
corner points are incident to a horizontal {\it and\/} a vertical
external link. 
We shall be interested here in allowed configurations in the FPL
model, in the sequel referred to as {\it FPL configurations},
with {\it periodic
boundary conditions}. These are FPL configurations where, around the
border of $Q_n$, every other external link of $Q_n$ 
is part of a polygon. The FPL configuration in
Figure~\ref{fpl24} is in fact a configuration with periodic boundary
conditions. 

\medskip

\medskip

It is well-known that FPL
configurations with periodic boundary conditions 
are in bijection with configurations in the six vertex
model with domain wall boundary conditions, which, in their turn, are
in bijection with alternating sign matrices (see, e.g., 
\cite[Sec.~3]{degier} for a description of these bijections).

Every FPL configuration with periodic boundary
condition defines a matching on the external links taken by the
polygons, by matching those which are on the same polygon. 
There has been a lot of interest recently in the enumeration of fully
packed loop configurations on $Q_n$ with periodic boundary
conditions, in which the matching on the external links is fixed.
To a big part this is due to the fact that it was 
(conjecturally) discovered that these numbers appear as the coordinates
of the groundstate vectors of certain Hamiltonians in the dense
$\text{O}(1)$ loop model. (See
\cite{degier} for a survey of these developments and conjectures). 

Although it is (probably) hopeless to expect a nice closed formula in
general, that is, for the number of FPL configurations with periodic
boundary conditions corresponding to a fixed matching, for an
arbitrary such matching, 
there exist several conjectures on these numbers for special matchings
(see, for example, \cite[Sec.~8]{wie}). In \cite{zuber}, Zuber added
several new ones, one of which he proved immediately in joint work with 
Di~Francesco and Zinn--Justin \cite{difra}. Another conjecture in this
direction for different boundary conditions, due to Mitra et al.\
\cite{Mitra}, 
was proved by de Gier in \cite[Sec.~5]{degier} (which is, in fact, the
inspiration for \cite{difra}, and also the present article). 
It is our purpose to prove two further conjectures from
\cite{zuber} in this paper, see Theorems~\ref{Z1} and \ref{Z2}. 

In all the proofs, the basic idea is to set up a bijection between the
FPL configurations in question and rhombus tilings of certain regions,
and then use known results on the enumeration of rhombus tilings to
conclude the proof. This is also the procedure which we shall follow
here. While in \cite{degier} the base of de Gier's result has been a theorem
of Ciucu and the second author \cite{ciuc}, and in \cite{difra} the
base of the result by Di~Francesco, Zinn--Justin and Zuber has been 
MacMahon's formula for plane partitions
contained inside a given box \cite{MM}, here it is
Stanley's hook-content formula \cite[Theorem~15.3]{StanAA} 
(see Theorem~\ref{ssyt} below) for the number of semistandard tableaux
of a given shape with bounded entries which is at the heart of our proofs. 
(We remark that this formula implies MacMahon's, see, e.g., \cite[proof
of Theorem~7.21.7]{StanBI}.)
In difference to \cite{difra,degier}, we are faced here with
an added difficulty in the proofs, as it is necessary to split the
enumeration problems considered here into several different
subcases. Another point worthy of note is the fact that our proofs use
Wieland's remarkable theorem of rotational symmetry \cite{wie} (see
Theorem~\ref{thm:wie} below) in an essential way, which is not
necessary in \cite{difra,degier}. That is, our proofs depend
crucially on the way the matching is ``placed around the square $Q_n$.''
We do in fact not know how to do the enumeration if we place the
matching in a different way around $Q_n$. On the other hand, it is
obvious that, using our approach, one can as well prove the
conjectures in Appendix~A of \cite{difra}, although we did not work
out the details. We do indeed hope that a refinement of the ideas
presented in this article will as well lead to a proof of
Conjectures~6 and 7 in \cite{zuber}. This is work currently in progress.

In the next section, we collect all notation and
the facts that we need in our
proofs. The proof of Conjecture~4 from \cite{difra} is then given in
Section~\ref{con4}, while the proof of Conjecture~5 from \cite{difra}
is the contents of Section~\ref{conj5}.

\section{Preliminaries\label{prel}}

We start by introducing the notation that we are going to use for
encoding FPL configurations and their associated matchings. The reader
should recall from the introduction that any FPL configuration defines
a matching on the external links taken by the
polygons, by matching those which are on the same polygon. We call
this matching the {\it matching associated to the FPL
configuration}. When we think of the matching as being fixed, and when
we consider all FPL configurations having this matching as associated
matching, we shall also speak of these FPL configurations as the ``FPL
configurations corresponding to this fixed matching."

We label the \( 4n \) external links around \( Q_{n} \)
in \( \mathbb Z/4n\mathbb Z \) clockwise starting from the left-most
link on the top side of the square, see Figure~\ref{fpl1}.
If \( A \) is an external link of the square, we denote by \( L(A) \)
its label and by \( LN(A) \) the representative of \( L(A) \) in
\( [-2n+1,2n] \). Throughout this paper, all the FPL configurations
that are considered are configurations which correspond to
matchings of either the even labelled external links or the odd
labelled external links.

Let \( M \) be any matching of the set of even (odd) labelled external
links. Let \( \tilde{M} \) be the ``rotated"
matching of the odd (even) external links defined by the property that
the links labelled $i$ and $j$ in $M$ are matched if and only if the
links labelled $i+1$ and $j+1$ are matched in $\tilde{M}$.
Let \( FPL(M) \) denote the number of FPL configurations
corresponding to the matching \( M \). Wieland \cite{wie}
proved the following surprising result. 

\begin{thm}[\sc Wieland]
\label{thm:wie}
For any matching \( M \) of the even (odd) labelled external links,
we have\[
FPL(M)=FPL(\tilde{M}).\]

\end{thm}

\begin{figure}[h]

{\centering \includegraphics{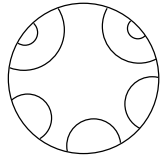} \par}

\caption{The chord diagram representation of a matching\label{fpl25}}
\end{figure}

In other terms, the number of FPL configurations corresponding to a
given matching is invariant under rotation of the ``positioning" of
the matching around the square. This being the case, we can represent
matchings in terms of chord diagrams of $2n$ points placed around a 
circle (see
Figure~\ref{fpl25} for the chord diagram representation of the
matching corresponding to the FPL configuration in
Figure~\ref{fpl24}). 

\medskip

In our proofs we shall also use the following observation of 
de Gier \cite[Lemma~8]{degier}. It is an assertion about the edges
that are taken by {\it any} FPL configuration if one makes certain
assumptions. 

\begin{lem}
\label{degier}Let $c$ be an FPL configuration which contains the edges
shown to the left of the implication symbol in Figure~{\em\ref{fpl16}}. We
assume furthermore that 
the top and the bottom edges do not belong to the same loop
and that one of the following two conditions is satisfied:

\begin{itemize} 
\item [i)]The middle edge belongs to a third loop.
\item [ii)]The middle edge is on the same loop as the top (bottom) 
one only if 
the loop contains the edge between
the left vertex of the top (bottom) 
edge and the right vertex of the middle edge.
\end{itemize}
Then $c$ contains all the edges on the right of the implication symbol. 

\begin{figure}[h]
\vspace{0.375cm}
{\centering \includegraphics{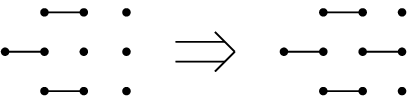} \par}
\caption{\label{fpl16}}
\end{figure}
\end{lem}
The following result is a consequence of an iterated use of 
Lemma~\ref{degier}. There, and in the sequel, when we speak of ``fixed
edges" we always mean edges that have to be occupied by {\it any} FPL
configuration under consideration.

\begin{lem}
\label{fixedg}Let \( A=A_{1},A_{2},\ldots ,A_{k}=B \) be
a sequence of external links, 
where $LN(A_i)=a+2i$ {\em mod} $4n$, for some fixed $a$, that is, the external
links $A_1,A_2,\dots,A_k$ comprise every second external link along the stretch
between $A$ and $B$ along the border of $Q_n$ (clockwise). 
Furthermore, we suppose that one of the following conditions holds:
\begin{enumerate}
\item \( A \) and \( B \) are both on the top side of \( Q_{n} \), that
is, \( 1\leq LN(A)<LN(B)\leq n \);
\item \( A \) is on the top side and \( B \) is on the right side of \( Q_{n} \),
that is, \( 1\leq LN(A)\leq n<LN(B) \) and \( n-LN(A)>LN(B)-(n+1) \);
\item \( A \) is on the left side and \( B \) is on the right side of
\( Q_{n} \), that is, \( n<LN(B)\leq 2n \) and \( -n<LN(A)\leq 0 \).
\begin{figure}
\psfrag{O}{\huge$O$}
\psfrag{A}{\huge$A$}
\psfrag{B}{\huge$B$}
\psfrag{A'}{\huge$A'$}
\psfrag{B'}{\huge$B'$}
{\centering \resizebox*{16cm}{!}{\includegraphics{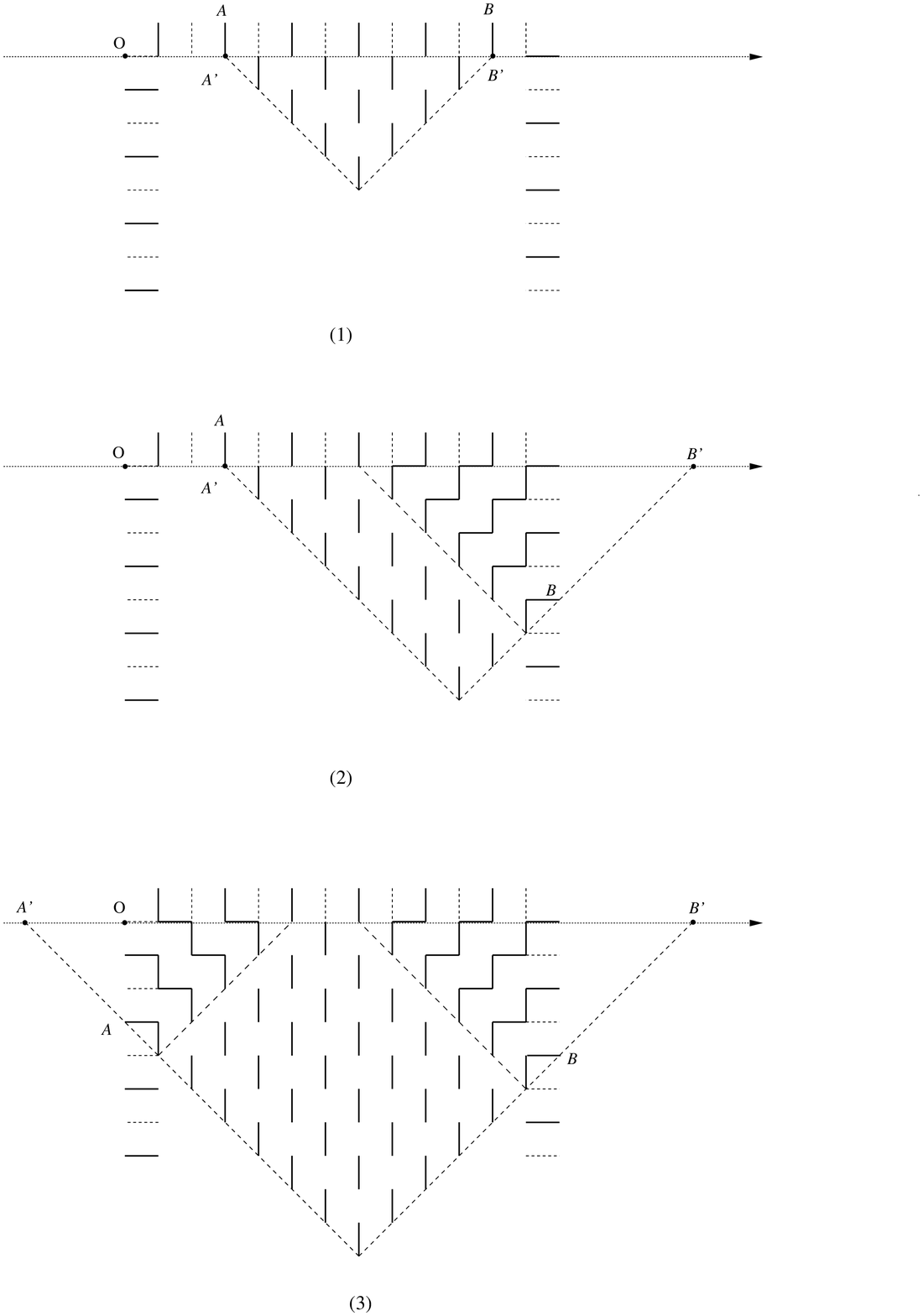}} \par}

\caption{The possible regions of fixed edges determined by a sequence
of external links
belonging to distinct loops\label{seq}}
\end{figure}

\end{enumerate}
For the FPL configurations for which the external links 
\( A_{1},A_{2},\ldots ,A_{k} \) belong to different loops,
the regions of fixed edges are {\em(}essentially{\em)} triangular
{\em(}see
Figure~{\em\ref{seq}} for illustrations; 
the {}``essentially'' refers to the fact that
in Cases~{\em(2)} and {\em(3)} parts of the triangle are cut
off{\em)}. More precisely, if one places
the origin \( O \) of the coordinate system one unit to the left of
the top-left corner of \( Q_{n} \), the coordinates of the triangle
are given in the following way: let \( A' \) and \( B' \) be the
points on the $x$-axis 
with \( x \)-coordinates \( LN(A) \) and \( LN(B) \), respectively,
then the region of fixed edges is given by the intersection of the square
\( Q_{n} \) and the {\em(}rectangular isosceles{\em)} 
triangle having the segment \( A'B' \) as basis. 

In Cases~{\em(2)} and {\em(3)}, 
the configurations are completely fixed as ``zig-zag" paths 
in the corner regions
of \( Q_{n} \) where a part of the triangle was cut off (see again
Figure~{\em\ref{seq}}). More precisely,
in Case~{\em(2)}, this region is the reflexion of the corresponding cut
off part of the triangle in the right side of \( Q_{n} \), and in
Case~{\em(3)} it is that region and also the reflexion of the corresponding
cut off part on the left in the left side of \( Q_{n} \). 
\end{lem}

\begin{figure}
{\centering \includegraphics{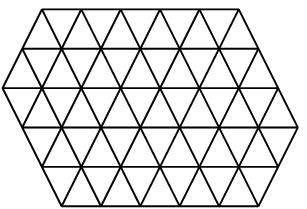} \par}

\caption{The hexagon $H(5,3,2)$\label{hexagon}}
\end{figure}

We next turn our attention to rhombus tilings of subregions of the
regular triangular lattice in the plane. 
Here, and in the sequel, by a
rhombus tiling we mean a tiling by rhombi of unit side lengths and
angles of $60^\circ$ and $120^\circ$. We first recall
MacMahon's theorem mentioned in the Introduction. 
Let \( H(p,q,r) \) be the hexagon with side lengths 
\( p ,  q ,  r ,  p ,  q ,  r \) (in clockwise order),  all of its
angles being $120^\circ$. We imagine $H(p,q,r)$ to be embedded in a
triangular lattice. See Figure~\ref{hexagon} for an illustration of
the hexagon \( H(5,3,2) \)). It is well known (see \cite{David-Tomei})
that rhombus tilings
of \( H(p,q,r) \) are in bijection with plane partitions contained
in a \( p\times q\times r \) box. 
The number of the latter plane
partitions was computed by MacMahon \cite[Sec.~429, $q \rightarrow 1$; 
proof in Sec.~494]{MM}.
Therefore we have the following theorem for the number \( h(p,q,r) \)
of rhombus tilings of \( H(p,q,r) \).

\begin{thm}[\sc MacMahon]
\label{plapar}Let \( p,q,r\in \mathbb N \). Then\[
h(p,q,r)=\frac{p\HH \, q\HH \, r\HH \, (p+q+r)\HH }{(p+q)\HH \, (p+r)\HH \, (q+r)\HH },\]
where \( n\HH:=(n-1)!\cdots 2!\,1! \) denotes
the \( n \)-th hyperfactorial.
\end{thm}

\begin{figure}
\psfrag{h1}{$h_1$}
\psfrag{h2}{$h_2$}
\psfrag{h3}{$h_3$}
\psfrag{h4}{$h_4$}
\psfrag{v1}{$v_1$}
\psfrag{v2}{$v_2$}
\psfrag{v3}{$v_3$}
\psfrag{v4}{$v_4$}

{\centering \includegraphics{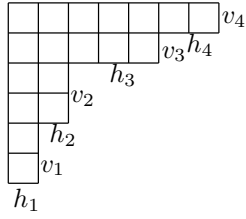} \par}

\caption{A Ferrers diagram\label{fpl18}}
\end{figure}

In the subsequent sections, we shall need a more general result for
regions which are indexed by partitions. We recall that a {\it partition} is
a vector $\la=(\la_1,\la_2,\dots,\la_\ell)$ of positive integers such
that $\la_1\ge\la_2\ge\dots\ge\la_\ell$. If there are repetitions
among the $\lambda_i$'s, then, for convenience, we shall sometimes use
exponential notation. For example, the partition $(3,3,3,2,1,1)$ will
also be denoted as $(3^3,2,1^2)$. To each partition $\la$, one
associates its {\it Ferrers diagram}, which is the left-justified
arrangement of cells with $\la_i$ cells in the $i$-th row,
$i=1,2,\dots,\ell$. See Figure~\ref{fpl18} for the Ferrers diagram of
the partition $(7,5,2,2,1,1)$. (At this point, the labels should be
disregarded.) The partition {\it conjugate to} $\la$
is the partition $\la'=(\la_1',\la_2',\dots,\la'_{\la_1})$, where
$\la_j'$ is the length of the $j$-th column of the Ferrers diagram of $\la$.
Given a partition $\la$, we write $(i,j)$ for the cell in the $i$-th
row and $j$-th column in the Ferrers diagram of $\la$,
\( 1\leq j\leq \lambda _{i} \). We use the notation 
\( u=(i,j)\in \lambda  \) to express the fact that $u$ is a cell of
(the Ferrers diagram of) $\la$. 
Given a cell $u$, we denote by \( c(u):=j-i \) the \emph{content}
of \( u \) and by \( h(u):=\lambda _{i}+\lambda '_{j}-i-j+1 \) the
\emph{hook length} of \( u \).

The enumeration result for rhombus tilings given in Theorem~\ref{thm:rhla}
below is a corollary of the hook-content formula for
semistandard tableaux of a given shape with bounded entries. Here,
{\it semistandard tableaux of shape $\la$}
are fillings of the cells of the Ferrers diagram of $\la$ with
positive integers such that the entries along rows are weakly
increasing and entries along columns are strictly increasing. See
Figure~\ref{fpl21} for a semistandard tableau of shape
$(7,5,2,2,1,1)$. 
We denote by \( SSYT(\lambda ,n) \) the set of semistandard tableaux
of shape \( \lambda  \) with entries less than or equal to \( n \).
Then Stanley's hook-content formula \cite[Theorem~15.3]{StanAA} reads
as follows.

\begin{thm}
\label{ssyt}Let \( \lambda  \) be a partition, and let 
\( n \) be a positive integer.
Then\[
|SSYT(\lambda ,n)|=\prod _{u\in \lambda }\frac{c(u)+n}{h(u)}.\]

\end{thm}
Given a partition \( \lambda  \), we are now going to define a region
in the regular triangular lattice which depends on $\la$. 
The bottom-right border of the
Ferrers diagram of $\la$ is a path consisting of positive unit horizontal and
vertical steps. This path determines the Ferrers diagram, and hence
the corresponding partition, uniquely. It may alternatively be described
by the sequence of lengths of its maximal horizontal and vertical
pieces \( (h_{1},v_{1},h_{2},v_{2},\ldots ,h_{k,},v_{k}) \). For
example this sequence for the partition in Figure~\ref{fpl18} is
\( (1,2,1,2,3,1,2,1) \).

\begin{figure}
\psfrag{h1}{$h_1$}
\psfrag{h2}{$h_2$}
\psfrag{v1}{$v_1$}
\psfrag{v2}{$v_2$}
\psfrag{hk}{$h_k$}
\psfrag{vk}{$v_k$}
\psfrag{r}{$r$}
\psfrag{l1}{$\lambda_1$}
\psfrag{l1'}{$\lambda_1'$}
\psfrag{r+}{$r+\sum_{i=1}^{k-1}v_i$}

{\centering \includegraphics{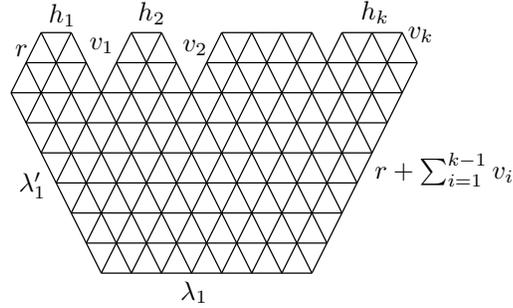} \par}

\caption{A hexagon with notches\label{fpl19}}
\end{figure}

Given a partition \( \lambda  \) and a nonnegative integer \( r \),
we define the region \( R(\lambda ,r) \) as a hexagon with some notches
along the top side. More precisely (the reader should consult
Figures~\ref{fpl18} and \ref{fpl19} in parallel, the latter showing
the region \( R(\lambda ,2) \), where
\( \lambda  \) is the partition in Figure~\ref{fpl18}), 
\( R(\lambda ,r) \) is the region
with base side \( \lambda _{1} \), bottom-left side \( \lambda _{1}' \),
top-left side \( r \), a top side with notches which will be explained
in just a moment, 
top-right side \( v_{k} \), and bottom-right side \( r+\sum _{i=1}^{k-1}v_{i} \).
Along the top side, we start with a horizontal piece of length \( h_{1} \),
followed by a notch of size \( v_{1} \), followed by a horizontal
piece of length \( h_{2} \), followed by a notch of size \( v_{2} \),
\( \ldots  \), and finally a horizontal piece of length \( h_{k} \).

We can now state the announced enumeration result for rhombus tilings
of the regions $R(\la,r)$. 

\begin{thm} \label{thm:rhla}
Given a partition \( \lambda  \) and a positive integer $r$, 
the number of rhombus tilings of
\( R(\lambda ,r) \) is given by \( |SSYT(\lambda ,r+\lambda _{1}')| \). 
\end{thm}
\begin{proof}
There is a standard bijection between our rhombus tilings and families
\( (P_{1},\ldots ,P_{\lambda _{1}}) \) of non-intersecting lattice
paths, where \( P_{i} \) is a path consisting of positive unit horizontal
and negative unit vertical steps from 
\( (i-\lambda _{i}',\lambda _{1}'-\lambda _{i}'+i+r) \)
to \( (i,i) \) (see, e.g., \cite{CEKZ,CK}). Here, ``non-intersecting" means 
the property that no two paths in a family have a point in common.
The bijection is obtained as follows. One places vertices in each of
the mid-points of edges along the base side of $R(\la,r)$, and as well
in each mid-point along the horizontal edges on the top of $R(\la,r)$.
The vertices of the top edges are subsequently connected to the
vertices along the base side by paths, by connecting the mid-points of
opposite horizontal edges in each rhombus of the tiling, see 
the tiling on the left of Figure~\ref{fpl20}.
Clearly, by construction, the paths are non-intersecting.
If the paths are slightly rotated, and deformed so that they become
rectangular paths, then one obtains families of paths with starting
and final points as described above. The family of paths which results
from our example rhombus tiling are shown on the right of
Figure~\ref{fpl20} (the labels should be ignored at the moment). 

\begin{figure}
\psfrag{1}{$1$}
\psfrag{2}{$2$}
\psfrag{3}{$3$}
\psfrag{4}{$4$}
\psfrag{5}{$5$}
\psfrag{6}{$6$}
\psfrag{7}{$7$}
\psfrag{8}{$8$}

{\centering \includegraphics{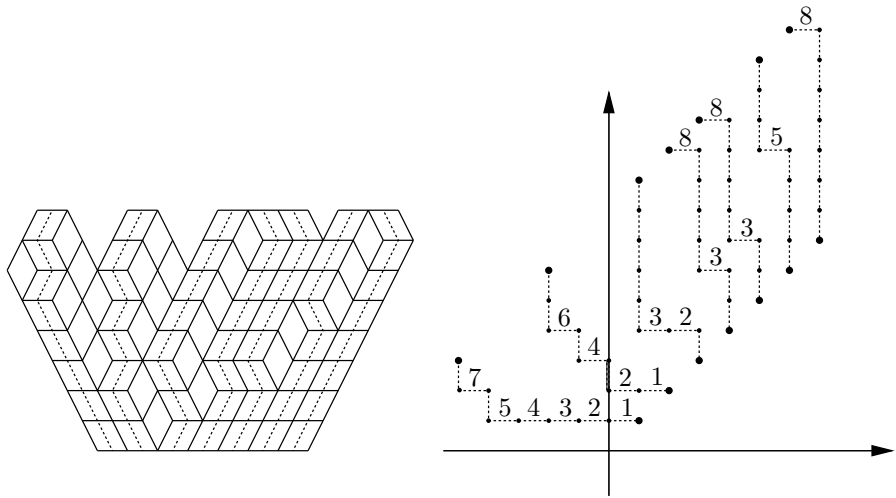} \par}

\caption{Bijection between rhombus tilings and non-intersecting
lattice paths\label{fpl20}}
\end{figure}
These families of non-intersecting lattice paths are, on the other
hand, in bijection with semistandard tableaux of shape \( \lambda  \)
with entries between \( 1 \) and \( r+\lambda _{1}' \) (see, e.g.,
\cite{FK}; it should be noted that Figure~8 there has to be reflected
in a vertical line to correspond to our picture). 
In this bijection, one labels the horizontal steps of the paths in
such a way, that a step from $(i,j)$ to $(i+1,j)$ gets the label
$j-i$, see Figure~\ref{fpl20}. A tableau is then formed by making the labels 
of the $j$-th path the entries of the $j$-th column of a tableau.
The tableau corresponding to the family of paths in 
Figure~\ref{fpl20} is shown in Figure~\ref{fpl21}.
\begin{figure}
\psfrag{1}{$1$}
\psfrag{2}{$2$}
\psfrag{3}{$3$}
\psfrag{4}{$4$}
\psfrag{5}{$5$}
\psfrag{6}{$6$}
\psfrag{7}{$7$}
\psfrag{8}{$8$}

{\centering \includegraphics{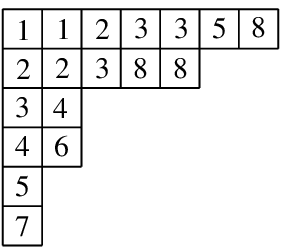} \par}

\caption{A semistandard tableau\label{fpl21}}
\end{figure}
\end{proof}
We need four special cases of this theorem in particular. We list
them explicitly in the following two lemmas. 
\begin{figure}
\psfrag{q}{\LARGE$q$}
\psfrag{r}{\LARGE\kern-5pt$r$}
\psfrag{p}{\LARGE$p$}
\psfrag{q-1}{\LARGE$q-1$}
\psfrag{r+1}{\LARGE$r+1$}
\psfrag{p+1}{\raise-2pt\hbox{\LARGE$p+1$}}

{\centering \resizebox*{6cm}{!}{\includegraphics{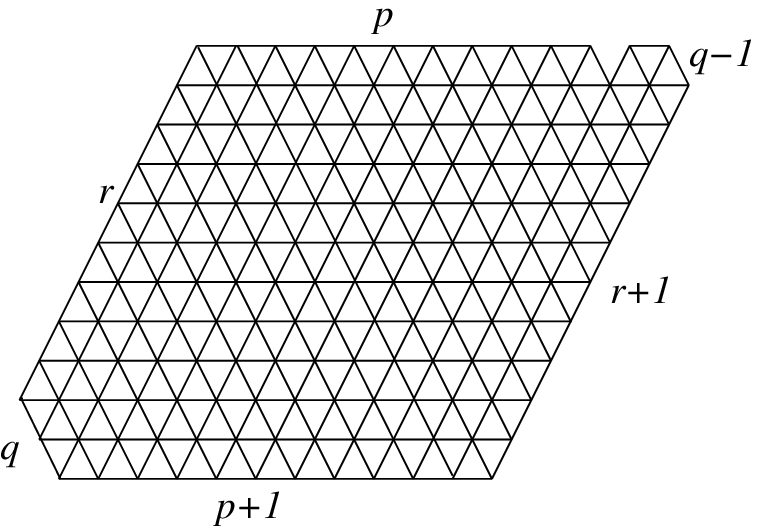}} \par}

\caption{\protect\( K(p,q,r)\protect \)\label{K(p,q,r)}}
\end{figure}

\begin{lem}
\label{k(p,q,r)}Let \( K(p,q,r) \) be the region \( R(((p+1)^{q-1},p),r) \)
{\em(}see Figure~{\em\ref{K(p,q,r)}} for an example with \( p=10 \), \( q=2 \)
and \( r=9 \){\em)}. Let \( k(p,q,r) \) be the number of rhombus tilings
of \( K(p,q,r) \). Then \[
k(p,q,r)=\frac{(p+q+r+1)\HH \, (p+1)\HH \, q\HH \, r\HH }{(p+q+2)\HH \, (p+r+2)\HH \, (q+r)\HH }(p+q)!\, (p+r)!\, q\, (p+1)\, (p+q+1).\]

\end{lem}
In order to state the next lemma more conveniently, we introduce the
following short notation:\[
a(p,q,r):=\frac{(p+q+r+2)\H \, r\H \, q\H \, (p+2)\H }{(p+q+4)\H \, (p+r+4)\H \, (q+r)\H }(p+q+1)!\, (p+q+2)!\, (p+r+3)!\, (p+r)!.\]

\begin{lem} \label{l(p,q,r)}
Let \( L(p,q,r)=R(((p+2)^{q},1,1),r-2) \), \(
M(p,q,r)=R(((p+2)^{q-1},p+1,1),r-1) \), 
and \( N(p,q,r)=R(((p+2)^{q-1},p),r) \) {\em(}see Figure~{\em\ref{L(p,q,r)}}
for illustrations of these regions and the corresponding Ferrers diagrams
in the case that \( p=4 \), \( q=3 \), and \( r=5 \){\em)}. 
Let \( l(p,q,r) \),
\( m(p,q,r) \), and \( n(p,q,r) \) be the numbers of rhombus tilings
of \( L(p,q,r) \), \( M(p,q,r) \), and \( N(p,q,r) \), respectively.
Then
\begin{enumerate}
\item \( l(p,q,r)=a(p,q,r)\frac{1}{2}(p+2)(p+3)(p+r+1)(p+r+2)r(r-1) \);
\item \( m(p,q,r)=a(p,q,r)(p+1)(p+3)(p+q+3)(p+r+1)qr \);
\item \( n(p,q,r)=a(p,q,r)\frac{1}{2}(p+1)(p+2)(p+q+2)(p+q+3)q(q+1) \).
\end{enumerate}

\begin{figure}
\psfrag{q+1}{\kern-5pt$q+1$}
\psfrag{r-1}{\kern-7pt$r-1$}
\psfrag{1}{$1$}
\psfrag{p}{$p$}
\psfrag{q-1}{\kern2pt$q-1$}
\psfrag{r+1}{$r+1$}
\psfrag{p+2}{$p+2$}
\psfrag{q}{\kern2pt$q$}
\psfrag{r}{$r$}
\psfrag{2}{$2$}
\psfrag{q+2}{\kern-6pt$q+2$}
\psfrag{r-2}{\kern-8pt$r-2$}
\psfrag{p+1}{$p+1$}
\psfrag{L(p,q,r)}{$L(p,q,r)$}
\psfrag{M(p,q,r)}{$M(p,q,r)$}
\psfrag{N(p,q,r)}{$N(p,q,r)$}

{\centering \includegraphics{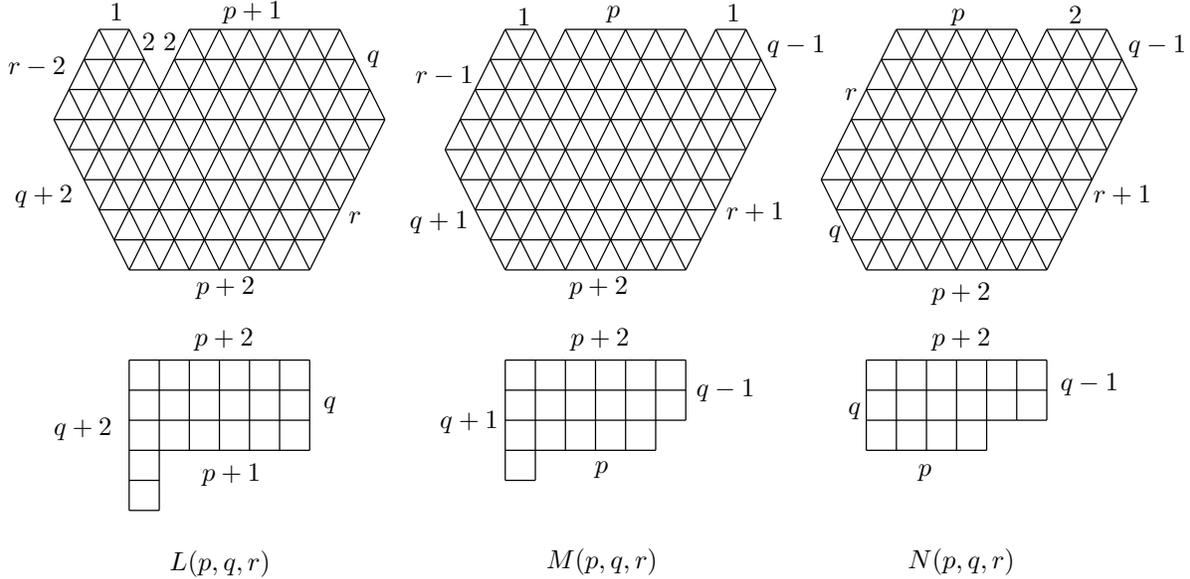} \par}

\caption{Some regions in the triangular lattice and the corresponding Ferrers
diagrams\label{L(p,q,r)}}
\end{figure}

\end{lem}

Finally, we need the following well-known result about local moves
applied to rhombus tilings, see, for example, 
\cite[Theorem~3.1]{SaToAA}.

\begin{lem}
\label{locmov}Let \( D \) be a simply connected region in the triangular
lattice. Then the local moves shown in Figure~{\em\ref{fig1}} act transitively
on the set of rhombus tilings of \( D \).
\begin{figure}[h]
{\centering \includegraphics{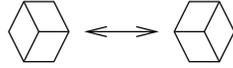} \par}

\caption{The local moves for rhombus tilings\label{fig1}}
\end{figure}
\end{lem}

\section{\label{con4}Proof of conjecture 4}

\begin{figure}
\psfrag{q}{\large\kern-2pt$q$}
\psfrag{r}{\large\kern1pt$r$}
\psfrag{p-1}{\large\kern5pt$p-1$}
\psfrag{Q}{\large\kern0pt$Q$}
\psfrag{R}{\large\raise2pt\hbox{\kern-1pt$R$}}
\psfrag{P}{\large\kern0pt$P$}

{\centering \includegraphics{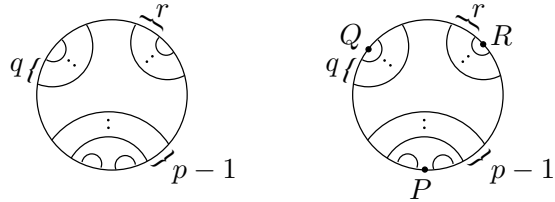} \par}

\caption{The matching in Theorem~\ref{Z1}\label{fpl9}}
\end{figure}

The goal of this section is to enumerate the fully packed loop
configurations whose associated matching is given by 
see the left half of Figure~\ref{fpl9}.

Note that we have \( n=p+q+r+1 \). We call 
the centres of the \( p-1 \), \( q \), and \( r \) nested arches
\( P \), \( Q \), and \( R \),
respectively, see the right half of Figure~\ref{fpl9}. 
We suppose first that \( p \) is even, and we let \( p=2s \).
Thanks to Wieland's Theorem~\ref{thm:wie}, we may place the linkage pattern
of the matching
arbitrarily around \( Q_{n} \). To prove the conjecture, we shall
make use of a particular placement, which we are going to explain
next.

We place the centre \( P \) on the external link labelled \( r+s+1 \).
This choice forces the other centres \( R \) and \( Q \) to be respectively
on the external links labelled by \( 5s+3r+4q+3=3n+(q-s) \) and \( 5s+2q+r+3=(2n+1)+(s-r) \).
In the sequel, by an FPL configuration we shall always mean an FPL
configuration corresponding to the particular matching in
Figure~\ref{fpl9}. Note that \( R \) is located on the bottom side
of the square if and only if \( s\geq q \), and similarly \( Q \)
is on the bottom side of the square if and only if \( s\geq r \).
Note also that we have a perfect right-to-left symmetry by exchanging
the roles of \( q \) and \( r \). 

A priori, we have three essentially distinct cases to deal with: 

\begin{enumerate}
\item \( R \) and \( Q \) are both on the bottom side; 
\item \( R \) is on the bottom side and \( Q \) is on the right side; 
\item \( R \) is on the left side and \( Q \) is on the right side. 
\end{enumerate}
We are going to concentrate on Case~(2). The Cases~(1) and (3) can
be treated in exactly the same way. In fact, all the claims that we
are going to make for Case~(2) are as well true for Cases~(1) and
(3).

\begin{figure}
\psfrag{O}{\raise3pt\hbox{\Large$O$}}
\psfrag{P}{\Large$P$}
\psfrag{Q}{\Large$Q$}
\psfrag{R}{\Large$R$}
\psfrag{R'}{\raise5pt\hbox{\Large\kern-4pt$R'$}}
\psfrag{R''}{\Large\kern-2pt$R''$}
\psfrag{Q}{\Large$Q$}
\psfrag{Q'}{\Large$Q'$}
\psfrag{Q''}{\Large$Q''$}
\psfrag{V}{\Large\kern-2pt$V$}
\psfrag{W}{\Large$W$}

{\centering \includegraphics{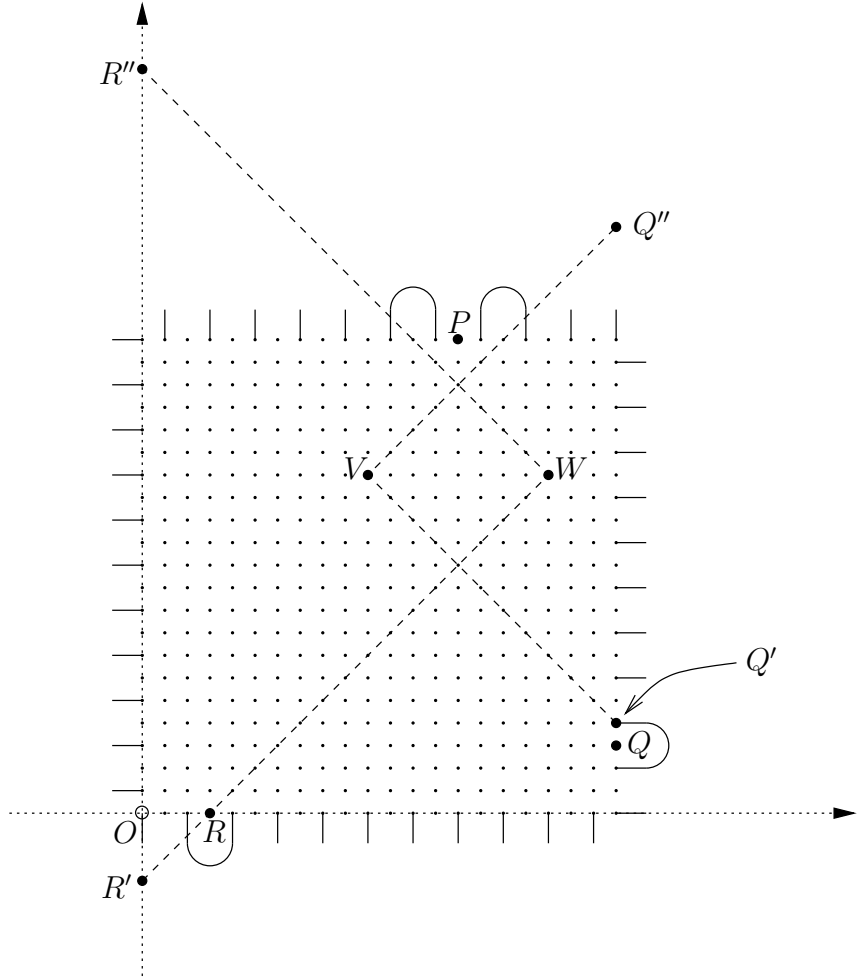} \par}

\caption{The region of fixed edges\label{VW}}
\end{figure}

\begin{figure}
\psfrag{P}{\Large$P$}
\psfrag{Q}{\Large$Q$}
\psfrag{R}{\Large$R$}

{\centering \includegraphics{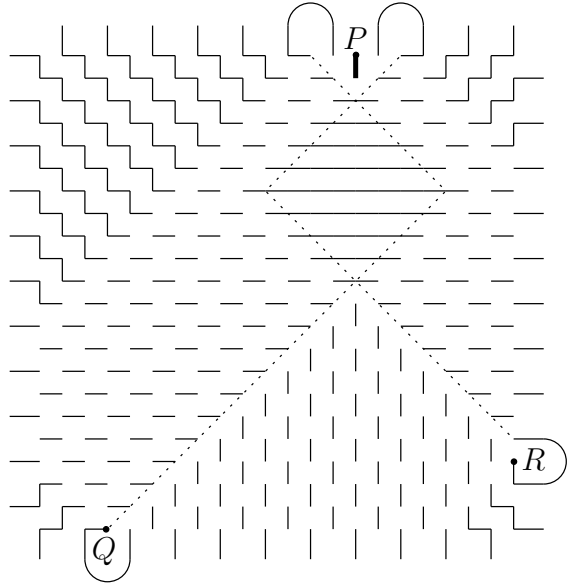} \par}

\caption{The set of fixed edges\label{edgfix}}
\end{figure}

For notational convenience we embed our square in the Cartesian plane
in such a way that the left side and the bottom side of the square
belong to the \( y \)-axis and the \( x \)-axis, respectively, with
unit length equal to the unit of the square, as shown in Figure~\ref{VW}.
We let \( R'':=(0,\, n-1+s+r+2) \), \( Q'':=(n-1,\, n-1+s+q+2) \),
\( R':=(0,q-s) \) and \( Q':=(n-1,r-s) \). Then, by Lemma~\ref{fixedg},
\( R' \) and \( R'' \) are two vertices of the triangle of fixed
edges determined by the distinct external links between \( P \) and \( R \),
and similarly for the triangle of fixed edges determined by the distinct
external links between \( P \) and \( Q \) (see Figure~\ref{VW}). The choice
of the position of the centres will ensure that the set of fixed edges
has certain useful properties.
In fact, if we call \( V \) and \( W \) the other vertices of these
triangles, as shown in Figure~\ref{VW}, we can easily see that they
have the same \( y \)-coordinate \( y_{V}=y_{W}=s+q+r-1 \) and that
\( x_{W}-x_{V}=p-2 \). So, the total set of fixed edges will be as
indicated in Figure~\ref{edgfix}.

If \( p \) is odd, we lose some of the symmetry of the case of even
\( p \) just discussed. Nevertheless, the argument remains essentially
valid. To be more precise, in this case 
we place \( P \) on the external link labelled \( r+s+1 \) where
\( s=\frac{1}{2}(p-1) \). The analogues of Figures~\ref{VW} and \ref{edgfix}
are essentially identical, except that the \( y \)-coordinates of
the vertices \( V \) and \( W \) differ by \( 1 \).

For both \( p \) even and odd, every vertex of the square belongs
to at least one fixed edge. If a vertex belongs to exactly \emph{one}
fixed edge we call it a \emph{free} vertex, and we say that two free
vertices are \emph{neighbours} if they can be joined by a non-fixed
edge. Now consider the vertical fixed edge just below \( P \) (marked
in bold-face in Figure~\ref{edgfix}). It is evident that the two
other edges emanating from its vertices have to be both on the right
or both on the left, otherwise we could not close the two small loops
next to \( P \). 

\begin{lem}
\label{ex}There exists a fully packed loop configuration for both
of these choices.
\end{lem}
\begin{proof}
This is very similar to the proof of \cite[Lemma~3]{difra} and is
hence omitted.
\end{proof}
\begin{thm}
\label{Z1}Let \( Z_{1}(p,q,r) \) be the number of fully packed loop
configurations determined by the matching in Figure~{\em\ref{fpl9}}. 
Then,\begin{multline*}
Z_{1}(p,q,r)  =  \frac{(p+q+r+1)\HH \, (p+1)\HH \, q\HH \, r\HH }
{(p+q+2)\HH \, (p+r+2)\HH \, (q+r)\HH }(p+q)!(p+r)! \\
 \times \Big ((p+1)(q(p+q+1)+r(p+r+1))+p(p+q+1)(p+r+1)\Big ).
\end{multline*}

\end{thm}
\begin{proof}
We treat the case where the edges emanating from the bold-face edge
are both on the left of it. This covers also the case where these
edges are both on the right, as is seen by an exchange of \( q \)
and \( r \). Then, following an idea of de Gier \cite[Sec.~5]{degier},
we draw a triangle around any \emph{free} vertex of our region in
such a way that two free vertices are neighbours if and only if the
corresponding triangles share an edge. This is illustrated in 
Figure~\ref{tria} for the example of Figure~\ref{edgfix} (the two added
fixed edges are marked in bold-face). 

\vspace{0.3cm}
{
\begin{figure}
\psfrag{A}{$A$}
\psfrag{B}{$B$}
\psfrag{C}{$C$}
\psfrag{D}{$D$}
\psfrag{E}{$E$}
\psfrag{F}{$F$}

{\centering \resizebox*{10cm}{10cm}{\includegraphics{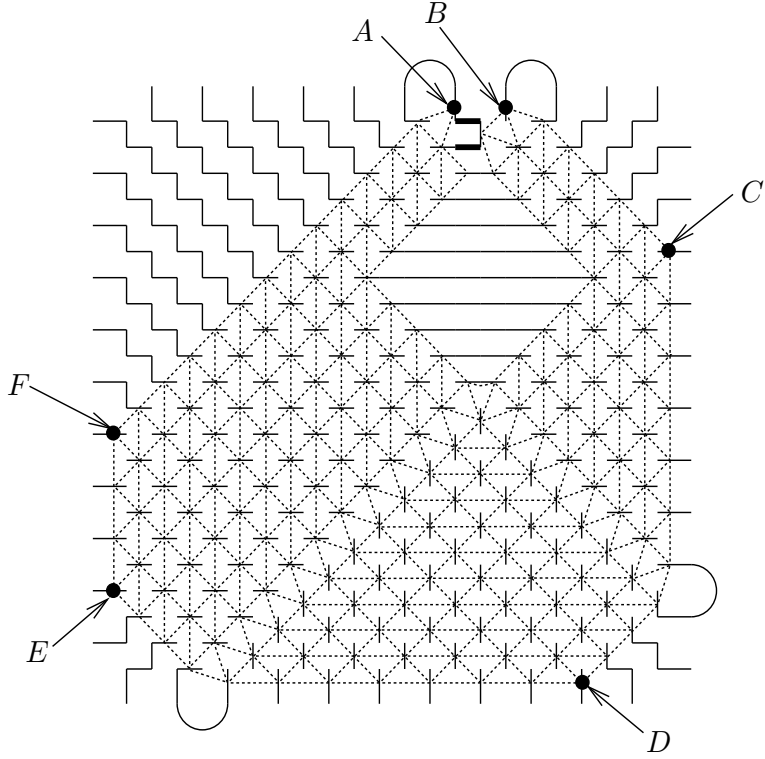}} \par}

\caption{The triangles around free vertices\label{tria}}
\end{figure}
\par}
\vspace{0.3cm}

Next we determine the dimensions of the region covered by these triangles.
Let \( A,B,C,D,E,F \) be the mid-points of the external links with labels
as indicated in the following table. See also Figure~\ref{tria}.
\\

{\centering \begin{tabular}{|c|c|}
\hline 
point&
link's label\\
\hline
\hline 
\( A \)&
\( s+r \)\\
\hline 
\( B \)&
\( s+r+2 \)\\
\hline 
\( C \)&
\( 3s+r+2q \)\\
\hline 
\( D \)&
\( 3s+3r+2q+2 \)\\
\hline 
\( E \)&
\( 7s+3r+2q+4 \)\\
\hline 
\( F \)&
\( 7s+3r+4q+6 \)\\
\hline
\end{tabular}\\
\par}

\bigskip

Given this information, the lengths of the sides of the region covered
by the triangles (measured in terms of triangle edges) are \( \overline{BC}=\frac{1}{2}((3s+r+2q)-(s+r+2))=s+q-1 \),
\( \overline{CD}=r+1 \), \( \overline{DE}=p+1 \), \( \overline{EF}=q+1 \)
and \( \overline{FA}=\frac{1}{2}((s+r)+4n-(7s+3r+4q+6))=s+r-1 \),
where in the last equality we have used the previously observed fact
that \( n=p+q+r+1=2s+q+r+1 \). Note the symmetry of these lengths
with respect to \( q \) and \( r \). 

After an appropriate deformation of the region in such a way that
it fits in a regular triangular lattice, we obtain a hexagon with
two {}``ears'', see Figure~\ref{ciao} for the result of the deformation
applied to the region of triangles in Figure~\ref{tria}.
\begin{figure}
\psfrag{q+1}{\kern-15pt\LARGE$q+1$}
\psfrag{s+r-1}{\kern-20pt\LARGE$s+r-1$}
\psfrag{s-1}{\LARGE$s-1$}
\psfrag{p-1}{\LARGE$p-1$}
\psfrag{s+q-1}{\LARGE$s+q-1$}
\psfrag{r+1}{\LARGE$r+1$}
\psfrag{p+1}{\LARGE$p+1$}

{\centering \resizebox*{6cm}{!}{\includegraphics{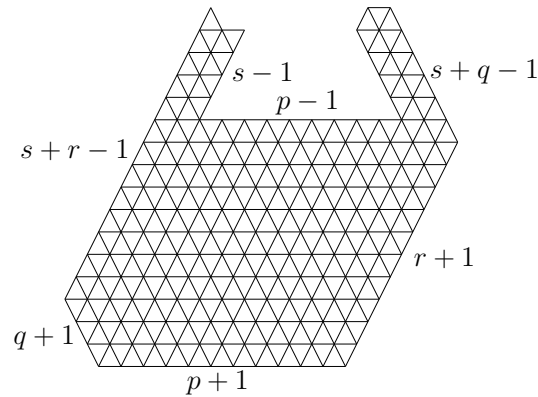}} \par}

\caption{The eared hexagon\label{ciao}}
\end{figure}

It is clear that any tiling of this region defines a fully packed
loop configuration just by drawing a segment between the two free
vertices corresponding to any tile. By Lemma~\ref{ex}, we know that
there is an FPL configuration for the case which we are discussing at the
moment. Any FPL configuration corresponds to a rhombus tiling of our
eared hexagon. By Lemma~\ref{locmov}, one can go from any rhombus
tiling of the eared hexagon to any other by the local moves 
shown in Figure~\ref{fig1}. It is easy to see that, under the
translation of FPL configurations into rhombus tilings, these moves
correspond to the local moves for FPL configurations shown in
Figure~\ref{local}. It is an important property of these latter moves that
they do not change the matching corresponding to the FPL
configurations. 
It follows that this correspondence establishes
a bijection between FPL configurations having the two prescribed edges
and rhombus tilings of the eared hexagon. 

\begin{figure}[h]

{\centering \resizebox*{6cm}{!}{\includegraphics{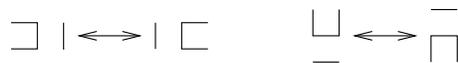}} \par}

\caption{Local moves for FPL configurations\label{local}}
\end{figure}
 
We now embark on the enumeration of the rhombus tilings of the hexagon
with two ears. In fact, in the left ear, and in a strip along the
left border, the tiles are uniquely determined (see Figure~\ref{fpl23}).
\begin{figure}
\psfrag{q+1}{\kern-15pt\LARGE$q+1$}
\psfrag{s+r-1}{\kern-20pt\LARGE$s+r-1$}
\psfrag{s-1}{\LARGE$s-1$}
\psfrag{p-1}{\LARGE$p-1$}
\psfrag{s+q-1}{\LARGE$s+q-1$}
\psfrag{r+1}{\LARGE$r+1$}
\psfrag{p+1}{\LARGE$p+1$}

{\centering \resizebox*{6cm}{!}{\includegraphics{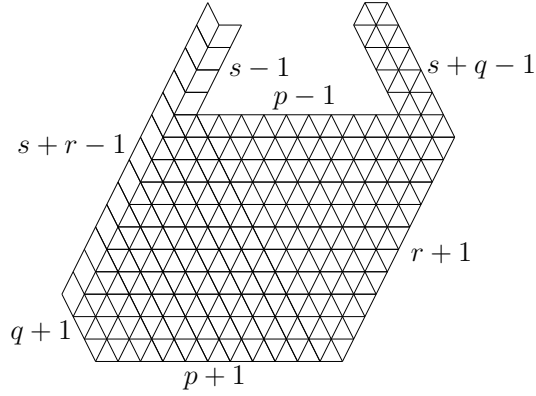}} \par}

\caption{Fixed rhombi in the eared hexagon\label{fpl23}}
\end{figure}
 Hence we can reduce our calculation to the number of tilings of a
hexagon with just one ear, see Figure~\ref{bon} (the shaded triangles
should be disregarded for the moment). We call this region \( F(p,q,r,s) \).
\begin{figure}
\psfrag{q}{\LARGE$q$}
\psfrag{r}{\kern-5pt\LARGE$r$}
\psfrag{p}{\LARGE$p$}
\psfrag{s-1}{\kern-17pt\LARGE$s-1$}
\psfrag{s+q-1}{\kern3pt\LARGE$s+q-1$}
\psfrag{r+1}{\LARGE$r+1$}
\psfrag{p+1}{\LARGE$p+1$}

{\centering \resizebox*{6cm}{!}{\includegraphics{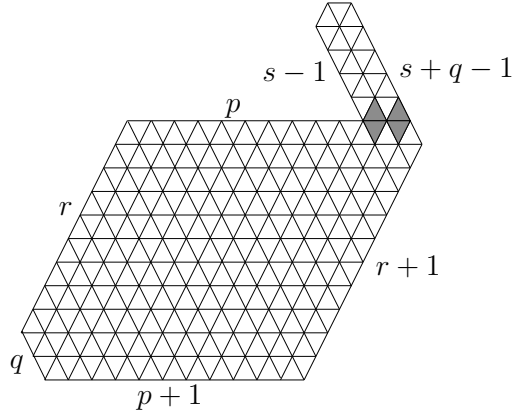}} \par}

\caption{Cutting the right ear\label{bon}}
\end{figure}

At this point, the crucial observation is that exactly one of the
two tiles which {}``connect'' the hexagon and the remaining ear
(these are the shaded tiles in Figure~\ref{bon}) must be chosen in
any tiling of \( F(p,q,r,s) \), otherwise our region would be cut
in two subregions with an odd number of triangles. If the right tile
is chosen, then our region is split into two separate hexagons, a
hexagon \( H(p+1,q,r) \) and a hexagon \( H(s-1,1,1) \) (see 
Section~\ref{prel} for the definition of \( H(p,q,r) \)). Hence the number
of rhombus tilings in this case is
\begin{align*}
h(p+1,q,r)&\,h(s-1,1,1)  =  \frac{(p+q+r+1)\H \, (p+1)\H \, q\H \, r\H }{(p+q+1)\H \, (p+r+1)\H \, (q+r)\H }s\\
 & =  \frac{(p+q+r+1)\H \, (p+1)\H \, q\H \, r\H }{(p+q+2)\H \, (p+r+2)\H \, (q+r)\H }(p+q)!\,(p+r)!\,(s(p+q+1)(p+r+1)).
\end{align*}

If the left tile is chosen, then the tiling of the ear is uniquely
determined, and the remaining region is \( K(p,q,r) \). Thanks to
Lemma~\ref{k(p,q,r)}, we already know the number of corresponding
rhombus tilings. Altogether, we obtain that the number of rhombus
tilings of \( F(p,q,r,s) \) is \[
f(p,q,r,s)=\frac{(p+q+r+1)\H\, (p+1)\H \, q\H \, r\H }{(p+q+2)\H \, (p+r+2)\H \, (q+r)\H }(p+q)!\,(p+r)!\,\Big (q(p+1)(p+q+1)+s(p+q+1)(p+r+1)\Big ).\]
 Hence, in total, there are \( Z_{1}(p,q,r)=f(p,q,r,p/2)+f(p,r,q,p/2) \)
FPL configurations. It is easy to verify that this agrees with our
claim. 

If \( p \) is odd, then, using the same approach, one shows that
in this case we have \( Z_{1}(p,q,r)=f(p,q,r,\frac{1}{2}(p-1))+f(p,r,q,\frac{1}{2}(p+1)) \),
which again agrees with our claim.
\end{proof}

\section{Proof of Conjecture 5\label{conj5}}

\begin{figure}
\psfrag{q}{$q$}
\psfrag{r}{$r$}
\psfrag{p-1}{$p-1$}
\psfrag{Q}{$Q$}
\psfrag{R}{\raise2pt\hbox{\kern-1pt$R$}}
\psfrag{P}{$P$}

{\centering \includegraphics{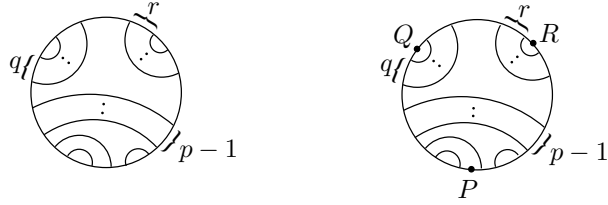} \par}

\caption{The matching in Theorem~\ref{Z2}\label{fant}}
\end{figure}

In this section we solve the problem of enumerating fully packed loop
configurations whose corresponding matching is described in
the left half of Figure~\ref{fant}. 
Our method of proof 
will be analogous to the one in the proof of Theorem~\ref{Z1},
though many more technicalities will occur.

\begin{figure}
\psfrag{O}{\LARGE$O$}
\psfrag{P}{\LARGE$P$}
\psfrag{Q}{\LARGE$Q$}
\psfrag{R}{\LARGE$R$}
\psfrag{R'}{\LARGE$R'$}
\psfrag{R''}{\kern-6pt\LARGE$R''$}
\psfrag{Q}{\LARGE$Q$}
\psfrag{Q'}{\LARGE$Q'$}
\psfrag{Q''}{\kern-4pt\LARGE$Q''$}
\psfrag{V}{\raise-5pt\hbox{\kern-8pt\LARGE$V$}}
\psfrag{W}{\raise-5pt\hbox{\kern3pt\LARGE$W$}}

{\centering \resizebox*{6cm}{!}{\includegraphics{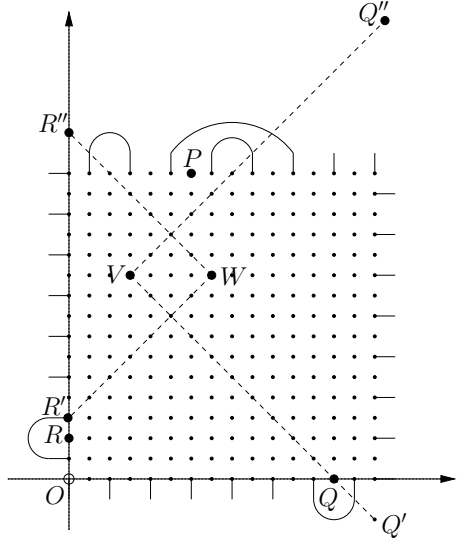}} \par}

\caption{The region of fixed edges\label{cha}}
\end{figure}

Again we first suppose \( p \) to be even, and we let \( p=2s \).
We call the centres of the \( p-1 \), \( q \), and \( r \) nested arches
\( P \), \( Q \), and \( R \),
respectively,  as before, see the right half of Figure~\ref{fant}. 
We place \( P \)
on the external link numbered \( r+s+2 \). Consequently \( Q \)
is on the external link labelled \( 5s+2q+r+6=(2n+1)+s-r+1 \), and \( R \)
is on the external link labelled \( 5s+4q+3r+6=3n+(q-s) \). We introduce the
same system of coordinates as in Section~\ref{con4}, and 
we let \( R':=(0,q-s) \),
\( R'':=(0,(n-1)+r+s-3) \), \( Q':=(n-1,r-s-1) \) and 
\( Q'':=(n-1,(n-1)+s+q-2) \).
Then these points determine the triangles of fixed edges (see 
Figure~\ref{cha}). Again, the vertices \( V \) and \( W \) of these triangles
have the same height, \( y_{V}=y_{W}=s+q+r-1 \), and 
we have \( x_{W}-x_{V}=p-2 \).

\begin{figure}
\psfrag{P}{\Large$P$}
\psfrag{Q}{\raise-1pt\hbox{\Large$Q$}}
\psfrag{R}{\Large\kern-1pt$R$}
{\centering \includegraphics{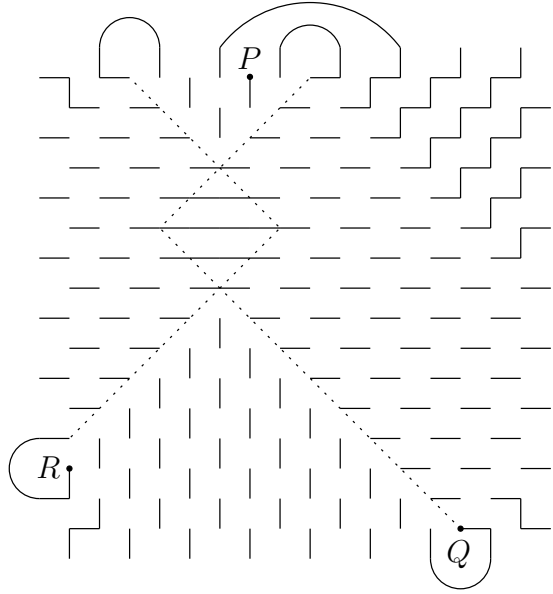} \par}

\caption{The set of fixed edges\label{fied}}
\end{figure}

At this point, by Lemma~\ref{fixedg}, we can easily draw the set
of fixed edges corresponding to our choice of the position of the
centres $P$, $Q$ and $R$ (see Figure~\ref{fied}).

\begin{figure}
\psfrag{P}{\Large$P$}
\psfrag{q+1}{\kern-6pt\large$q+1$}
\psfrag{r-1}{\large$r-1$}
\psfrag{1}{\kern-1pt\large$1$}
\psfrag{p}{\large$p$}
\psfrag{q-1}{\large$q-1$}
\psfrag{r+1}{\large$r+1$}
\psfrag{p+2}{\large$p+2$}
\psfrag{q}{\kern-1pt\large$q$}
\psfrag{r}{\large$r$}
\psfrag{2}{\large$2$}
\psfrag{q+2}{\kern-5pt\large$q+2$}
\psfrag{r-2}{\large$r-2$}
\psfrag{p+1}{\large$p+1$}
\psfrag{r+s-2}{\kern-14pt\large$r+s-2$}
\psfrag{s-2}{\kern-4pt\large$s-2$}
\psfrag{q+s}{\large$q+s$}
\psfrag{s-1}{\kern-9pt\large$s-1$}
\psfrag{q+s-1}{\large$q+s-1$}
\psfrag{s}{\kern-1pt\large$s$}

{\centering \includegraphics{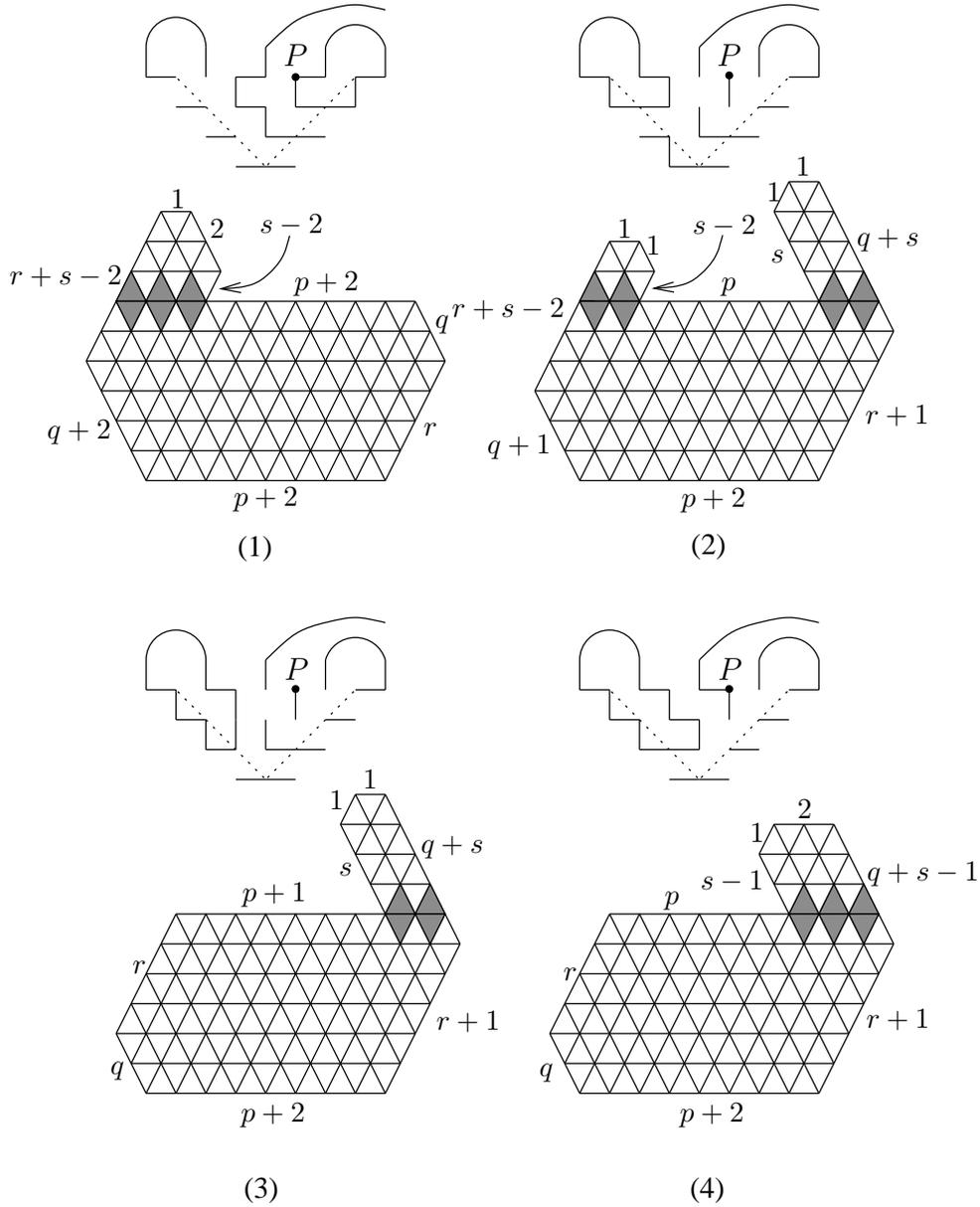} \par}

\caption{The possible local configurations and the corresponding regions\label{config}}
\end{figure}

Again, any vertex of the square grid belongs to at least one fixed
edge. However, here we have to split our problem into four cases,
according to the local configuration around the centre \( P \). There
are indeed four mutually exclusive possible configurations near \( P \),
see Figure~\ref{config}.

\begin{lem}
\label{ex2}There exists a fully packed loop configuration for any
of the four local configurations shown in Figure~{\em\ref{config}}.
\end{lem}
\begin{proof}
Again the proof is omitted since it is very similar to the proof of
\cite[Lemma~3]{difra}.
\end{proof}
\begin{thm}
\label{Z2}
Let \( Z_{2}(p,q,r) \) be the number of fully packed loop configurations
determined by the matching in Figure~{\em\ref{fant}}. Then
\begin{multline*}
Z_{2}(p,q,r)  =  \frac{a(p,q,r)}{2}(p+2)\Big ((p+1)(p+q+3)(p+r+1)(p(p+r+2)+pq+4qr)+\\
  +2p(p+q+3)q(p+r+1)(p+q+2)+2(p+1)(p+q+3)(p+r+1)(p+r+2)r+\\
  +(p+3)(p+r+1)(p+r+2)r(r-1)+(p+1)(p+q+2)(p+q+3)q(q+1)\Big ).
\end{multline*}
 
\end{thm}
\begin{proof}
Again, we only discuss the case where \( p \) is even, since the other
case is completely analogous. We already know that \( Z_{2}(p,q,r) \)
is the sum of four parts corresponding to the local configurations
described in Figure~\ref{config}. In each of these cases we proceed
similarly as in the proof of Theorem~\ref{Z1}. We draw a triangle
around any free vertex in such a way that triangles share an edge
if and only if the corresponding vertices are adjacent. After a suitable
deformation of the corresponding regions, we can fit these regions
in the regular triangular lattice. Again, there are certain parts which are
covered by forced tiles and which may hence be eliminated. See 
Figure~\ref{config} for an illustration of the resulting regions in each
of the four cases (the shaded triangles should be disregarded for
the moment). We now observe that Lemmas~\ref{locmov} and \ref{ex2}
imply that in each of the four cases there is a bijection between
the fully packed loop configurations and rhombus tilings of the corresponding
region in the triangular lattice. So the number of fully packed loop
configurations is equal to the sum of the number of tilings of these
four regions. We treat them one at a time. Below, the numbers (1)--(4)
refer to the corresponding numbers in Figure~\ref{config}.
\begin{enumerate}
\item Consider the three shaded tiles. It is easy to see that exactly two
of them should be chosen in any tiling of this region. If we choose
the first two, then the region is split into two hexagons, a hexagon
\( H(s-2,1,2) \) and a hexagon \( H(p+2,q,r) \). If we choose the
first and the third, then the region is split into a hexagon \( H(s-1,1,1) \)
and a region \( K(p+1,r,q) \). Finally, if we choose the second and
the third tile, our region reduces to a region \( L(p,q,r) \). Hence
the total number of rhombus tilings in this case is\[
\frac{1}{2}s(s-1)\cdot h(p+2,q,r)+s\cdot k(p+1,r,q)+l(p,q,r).\]

\item Any tiling of the second region must contain exactly one of the two
shaded tiles on the left and exactly one of the two shaded tiles on
the right. If we choose the first on both the right and the left,
then we obtain a hexagon \( H(s-2,1,1) \) and a region \( K(p+1,q,r) \).
If we choose the first on the left and the second on the right, we
split the region into three hexagons, a hexagon \( H(s-2,1,1) \),
a hexagon \( H(p+2,q,r) \), and a hexagon \( H(s,1,1) \). If we choose
the second on the left and the first on the right, we obtain a region
\( M(p,q,r) \). Finally, if we choose the second on both the left
and the right, we split the region into a hexagon \( H(s,1,1) \) and
a region \( K(p+1,r,q) \). Hence the number of rhombus tilings in
this case is\[
(s-1)\cdot k(p+1,q,r)+(s^{2}-1)\cdot h(p+2,q,r)+m(p,q,r)+(s+1)\cdot k(p+1,r,q).\]

\item Again, we have to choose exactly one of the two shaded tiles. If we
choose the first one, we get a region \( K(p+1,q,r) \), and if we
choose the second one, we get two disjoint hexagons, a hexagon \( H(s,1,1) \)
and a hexagon \( H(p+2,q,r) \). Hence the number of rhombus tilings
in this case is\[
k(p+1,q,r)+(s+1)\cdot h(p+2,q,r).\]

\item In this last case we have to choose exactly one of the three tiles.
If we choose the first one, we get a region \( N(p,q,r) \). If we
choose the second, we obtain a hexagon \( H(s-1,1,1) \) and a region
\( K(p+1,q,r) \), and if we choose the third one, we get a hexagon
\( H(s-1,1,2) \) and a hexagon \( H(p+2,q,r) \). Hence the number
of rhombus tilings in this case is \[
n(p,q,r)+s\cdot k(p+1,q,r)+\frac{1}{2}s(s+1)\cdot h(p+2,q,r).\]

\end{enumerate}
Putting everything together, the total number \( Z_{2}(p,q,r) \)
is given by the sum of these four expressions, that is
\begin{multline*}
Z_{2}(p,q,r)=(2s^{2}+s)\cdot h(p+2,q,r)+2s\cdot k(p+1,q,r)+
(2s+1)\cdot k(p+1,r,q)\\
+l(p,q,r)+m(p,q,r)+n(p,q,r).
\end{multline*}
Now note that, by Theorem~\ref{plapar} and Lemma~\ref{k(p,q,r)}, there
hold
\begin{align*}
h(p+2,q,r) & =  a(p,q,r)(p+q+2)(p+q+3)(p+r+1)(p+r+2),\\
k(p+1,q,r) & =  a(p,q,r)(p+2)(p+q+3)q(p+r+1)(p+q+2),
\end{align*}
and\[
k(p+1,r,q)=a(p,q,r)(p+2)(p+q+3)(p+r+1)(p+r+2)r.\kern13pt\]
If we also recall Lemma~\ref{l(p,q,r)}, then it follows that
\begin{align*}
Z_{2}&(p,q,r)  =  a(p,q,r)\Big (\frac{p}{2}(p+1)(p+q+2)(p+q+3)(p+r+1)(p+r+2)\\
 &   \kern3.5cm +p(p+2)(p+q+3)q(p+r+1)(p+q+2)\\
 &   \kern3.5cm +(p+1)(p+2)(p+q+3)(p+r+1)(p+r+2)r\\
 &   \kern3.5cm +\frac{1}{2}(p+2)(p+3)(p+r+1)(p+r+2)r(r-1)\\
 &   \kern3.5cm +(p+1)(p+3)(p+q+3)(p+r+1)qr\\
 &   \kern3.5cm +\frac{1}{2}(p+1)(p+2)(p+q+2)(p+q+3)q(q+1)\Big )\\
 & =  \frac{a(p,q,r)}{2}(p+2)\Big ((p+1)(p+q+3)(p+r+1)(p(p+r+2)+pq+4qr)\\
 &   \kern1cm +2p(p+q+3)q(p+r+1)(p+q+2)+2(p+1)(p+q+3)(p+r+1)(p+r+2)r\\
 &   \kern1cm +(p+3)(p+r+1)(p+r+2)r(r-1)+(p+1)(p+q+2)(p+q+3)q(q+1)\Big ),
\end{align*}
which agrees with the expression in the assertion of the theorem.
\end{proof}

\section*{Acknowledgement}
We thank Jean-Bernard Zuber for several corrections and useful remarks
on an earlier version of this paper.

\end{document}